# A class of new stable, explicit methods to solve the non-stationary heat equation

Endre Kovács

University of Miskolc, Institute of Physics and Electrical Engineering,

kendre01@gmail.com

Highlights

- New first and second order schemes to solve the heat or diffusion equation.
- Unconditional stability for the heat equation, no stepsize requirements.
- Explicit and one-step methods, easy to implement and parallelize.
- Handy to apply regardless of space dimensions and grid irregularity.

**Abstract:** We present a class of new explicit and stable numerical algorithms to solve the spatially discretized linear heat or diffusion equation. After discretizing the space and the time variables like conventional finite difference methods, we do not approximate the time derivatives by finite differences, but use constant-neighbour and linear neighbour approximations to decouple the ordinary differential equations and solve them analytically. During this process, the timestep-size appears not in polynomial, but in exponential form with negative exponents, which guarantees stability. We compare the performance of the new methods with analytical and numerical solutions. According to our results, the methods are first and second order in time and can be much faster than the commonly used explicit or implicit methods, especially in the case of extremely large stiff systems.

## 1. Introduction and the studied problem

It is well known that the simplest diffusion and Fourier-type heat conduction phenomena are described by a second-order linear parabolic partial differential equation (PDE), the so-called heat equation:

$$\frac{\partial u}{\partial t} = \alpha \nabla^2 u + q, \tag{1}$$

where $u$ is the temperature (the concentration in the diffusion-equation), $\alpha = k/(c\rho) > 0$ is the thermal diffusivity, $q$, $k$, $c$, and $\rho$ is the intensity of heat sources (radiation, chemical reactions radioactive decay, etc.), heat conductivity, specific heat and (mass) density, respectively.

It is well known that this equation can describe the conduction of heat in many different systems from power plants to building walls [1]. Its generalizations like the convection (or advection)-diffusion equation and reaction-diffusion equation can describe the diffusion of particles in chemical and biological systems as well as in electronic devices [2]. Moreover, its analogs are widely used to treat fluid flow through porous media [3], like ground water or crude oil and natural gas in reservoirs.

To solve this equation numerically, the most typical starting step is the same as in the standard *method of lines*. The most common finite difference scheme to discretize the second space derivatives [4, p. 48] is the second order central difference formula:

$$\frac{\partial^2}{\partial x^2} f(x_i, t_j) \approx \frac{\frac{f(x_{i+1}, t_j) - f(x_i, t_j)}{\Delta x} + \frac{f(x_{i-1}, t_j) - f(x_i, t_j)}{\Delta x}}{\Delta x}.$$

For a realistic system, the discretization must reflect the material properties of the system, thus we perform it more generally than standard numerical analysis textbooks usually do, e.g. we don't consider $\alpha$, $k$, $c$ and $\rho$ as spatially uniform. Using the above formula for a one dimensional, equidistant grid, we obtain

$$\left.\frac{\partial u}{\partial t}\right|_x = \frac{1}{\Delta x \cdot c|_x \rho|_x}\left(k|_{x+\frac{\Delta x}{2}} \cdot \frac{u(x+\Delta x)-u(x)}{\Delta x} + k|_{x-\frac{\Delta x}{2}} \cdot \frac{u(x-\Delta x)-u(x)}{\Delta x}\right) + q|_x.$$

Now we change to cell variables:

$$\frac{du_i}{dt} = \frac{A}{c_i \rho_i A \Delta x}\left(k_{i,i+1} \cdot \frac{u_{i+1}-u_i}{\Delta x} + k_{i-1,i} \cdot \frac{u_{i-1}-u_i}{\Delta x}\right) + Q_i,$$

where $u_i$ is the temperature of the cell $i$, $C = c \cdot m = c\rho V$ is the heat capacity of the cell in $[J/K]$ units ($m$ is the mass, $V = A\Delta x$ is the volume of the cell). We introduce two other quantities, the heat source term $Q$,

$$Q_i = \frac{1}{V_i}\int_{V_i} q\,dV \approx q, \text{ in } \left[\frac{K}{s}\right] \text{units},$$

and the thermal resistance $R_{ij} = \frac{\Delta x}{k_{ij}A}$ in $[K/W]$ units. In case of an irregular grid, the distances between the cell-centres are $d_{ij} = (\Delta x_i + \Delta x_j)/2$ and the resistances can be approximated as $R_{ij} \approx \frac{d_{ij}}{k_{ij}A_{ij}}$.

Using these notations we obtain

$$\frac{du_i}{dt} = \frac{u_{i-1}-u_i}{R_{i-1,i}C_i} + \frac{u_{i+1}-u_i}{R_{i+1,i}C_i} + Q_i.$$

For a standard one dimensional homogeneous system with an equidistant grid, this would take a more familiar form:

$$\frac{du_i}{dt} = \alpha\frac{u_{i-1}-2u_i+u_{i+1}}{(\Delta x)^2} + Q_i.$$

We prefer to use the ODE system for a general (possibly unstructured) grid, which gives the time derivative of each temperature independently of any coordinate-system:

$$\frac{du_i}{dt} = \sum_{j\neq i}\frac{u_j-u_i}{R_{i,j}C_i} + Q_i.$$

In a matrix-form, we have

$$\frac{d\vec{u}}{dt} = M\vec{u} + \vec{Q}, \qquad (2)$$

where an off-diagonal $m_{ij} = 1/(R_{ij}C_i)$ element of the $M$ matrix can be nonzero only if the cells $i$ and $j$ are neighbours. From this point, all summations are going over the neighbours of the actual cell, which will be denoted by $j \in n(i)$. Unless stated otherwise, we consider closed (zero Neumann) boundary conditions, i.e. the edge of the examined domain is thermally isolated regarding conductive type heat transfer. In order to help the reader to imagine, we present the arrangement of the variables in Fig 1. for a 2D system of 4 cells. We emphasize that the shape and arrangement of the cells are not necessarily regular.

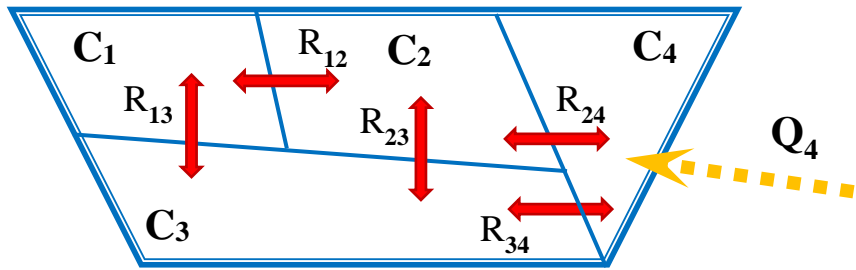

Figure 1 Notations in the case of 4 cells. The outer double line represents thermal isolation and the inner thin lines represent cell borders. The red double arrows are for conductive heat transport through the resistances $R_{ij}$ while $Q_4$ is an external heat source.

For this simple system, the differential equation in matrix form is the following:

$$\frac{d}{dt}\begin{pmatrix} u_1 \\ u_2 \\ u_3 \\ u_4 \end{pmatrix} = \begin{pmatrix} \frac{-1}{C_1 R_{12}} + \frac{-1}{C_1 R_{13}} & \frac{1}{C_1 R_{12}} & \frac{1}{C_1 R_{13}} & 0 \\ \frac{1}{C_2 R_{12}} & \frac{-1}{C_2 R_{12}} + \frac{-1}{C_2 R_{23}} + \frac{-1}{C_2 R_{24}} & \frac{1}{C_2 R_{23}} & \frac{1}{C_2 R_{24}} \\ \frac{1}{C_3 R_{13}} & \frac{1}{C_3 R_{23}} & \frac{-1}{C_3 R_{13}} + \frac{-1}{C_3 R_{23}} + \frac{-1}{C_3 R_{34}} & \frac{1}{C_3 R_{34}} \\ 0 & \frac{1}{C_4 R_{24}} & \frac{1}{C_4 R_3} & \frac{-1}{C_4 R_{24}} + \frac{-1}{C_4 R_{34}} \end{pmatrix} \begin{pmatrix} u_1 \\ u_2 \\ u_3 \\ u_4 \end{pmatrix} + \begin{pmatrix} Q_1 \\ Q_2 \\ Q_3 \\ Q_4 \end{pmatrix}.$$

One can see that the size of the matrix grows quadratically with the number of cells, thus the number of elements of the matrix is inversely proportional with the 4th power of the diameter of the blocks (for a fixed 2D system) and with the 6th power of the diameter in 3D. On the other hand, in real systems, the physical properties like the specific heat and the heat conductivity can widely vary from point to point [1, p. 15]. Therefore the magnitude of the matrix elements and thus the eigenvalues may have a range of several orders of magnitude, which means that the problem can be severely stiff. As the conventional explicit methods are only conditionally stable, they are inappropriate to solve these kinds of problems, because unacceptably small timesteps would be required. That is why implicit methods are almost exclusively used, as they have good stability properties. However, they require the solution of an algebraic equation system at each time-step, usually by iterative methods. As the parallelization of these methods is nontrivial, the solution is still rather slow, especially when the matrix is huge, and in more than one space dimension, where the matrix is not tridiagonal. The simulation is especially challenging in the oil industry, where they have to deal with multiphase flow in underground reservoirs with a diameter as large as tens of kilometres containing heterogeneous fractured porous media [5]. To overcome the difficulties, sophisticated techniques have been developed, among which we mention only a few. The most time-consuming part of the simulation is the solution of the sequences of large, often ill-conditioned linear systems of equations. To accelerate the convergence rate and to avoid divergence, more and more effective preconditioners are applied, e.g. the two stage preconditioner of White et al., see [6] and the references therein. Another powerful idea is the application of a multigrid system, i.e. a hierarchy of discretization. First, geometric multigrid methods appeared, which used grids with different resolutions when the variables of the PDE were discretized. On the other hand, algebraic multigrid (AMG) methods construct their hierarchy of operators directly from the system matrix thus they can be used regardless of the geometry of the simulated system. Nowadays especially AMG methods are used by some commercial and open-source software as well [7], but hybrid solutions have also been developed and applied, which combine algebraic and geometric methods [8], [9]. Moreover, dynamic local grid refinement methods, like Algebraic Dynamic Multilevel (ADM) method, employ high-resolution grids only where and when necessary [10].

Not everyone is aware that there are explicit methods with better stability properties. The textbook example is the Dufort-Frankel scheme, but we have to mention Runge–Kutta–Chebyshev methods [11], the Hopscotch method [12], [13], Alternating Direction Explicit (ADE) and Alternating Group Explicit (AGE) methods [14]–[16], the D'Yakonov fully explicit variant of the iterative alternating decomposition method [17] and the positivity preserving schemes of other research groups [18], [19]. However, the investigation and application of these methods are relatively rare (see for example [20]–[23]) and implicit methods are widely considered to be superior [24][25]. One of the main reasons of this underestimation must be that these explicit methods have other disadvantages. For example, they can only be conditionally convergent or consistent, less accurate or complicated to code, or they can hardly be applied for irregular grids [26], [27, p. 224], [28]–[30].

One may conclude that known methods are not really convenient for this purpose. In order to fill this gap, we started to elaborate a class of fundamentally new explicit methods. The simplest, first order member of this class has been already published [31]. In the next section, we introduce the new algorithms and describe their properties. We have been comparing the results produced by our methods with a few analytical and plenty of numerical results as well. We present here one analytical case in section 3 and two numerical cases in sections 4 and 5. Finally, we summarize the conclusions and the possible directions of our subsequent research in section 6.

## 2. The description of the methods

Now we introduce the new methods to solve the ODE system (2) through the following steps.

**1a)** At the predictor step, we make an approximation: when we calculate the new value of a variable $u_i^{n+1}$, we neglect that other variables are also changing during the timestep $h = \Delta t$. It means that we consider $u_j$ a constant if $j \neq i$, so we have to solve *uncoupled*, linear ODEs:

$$\frac{du_i}{dt} = a_i - \frac{u_i}{\tau_i}, \tag{3}$$

where we introduced $a_i = \sum_{j \in n(i)} m_{ij} u_j^n + Q_i$, and the characteristic time or time-constant of the cell $\tau_i = -1/m_{ii} = C_i / \sum_{j \in n(i)} 1/R_{ij}$, while $u_j^n$ are the temperatures at the beginning of the *n*th timestep. Although this approximation may look similar to the explicit Euler method, the crucial difference is that we do not fix the actual variable $u_i$ to its initial value (we fix only the neighbours), but keep it as a variable. It means that we still have *differential* equations and not *difference* equations.

**1b)** It is easy to solve (3) **analytically**. We use this solution at the end of the timestep as the predictor values:

$$u_i^{n+1,\text{pred}} = u_i^n \cdot e^{-\frac{h}{\tau_i}} + a_i \tau_i \cdot \left(1 - e^{-\frac{h}{\tau_i}}\right). \tag{4}$$

If we accepted these predictor values as final results at the end of the timestep, we would obtain a one-stage first order method, which could be called "**constant-neighbour** method" or CN1, while the concrete values can be denoted as $u_i^{\text{CN1}}$ or simply $u_i^{\text{CN}}$.

**2a)** During the corrector step, we make a more realistic assumption that the neighbouring variables $u_j$, of the actual cell $u_i$ are changing **linearly**. For each cell $i$ we can introduce an aggregated or effective slope for the neighbours:

$$s_i = \frac{a_i^{\text{pred}} - a_i}{h},$$

where $a_i^{\text{pred}} = \sum_{j \in n(i)} m_{ij} u_j^{n+1,\text{pred}}$ contains the predictor values of all the neighbours. Using this approximation we obtain a new uncoupled ODE system:

$$\frac{du_i}{dt} = s_i t + a_i - \frac{u_i}{\tau_i} \tag{5}$$

**2b)** It is still not difficult to solve eq. (4) analytically. We use the solution of it at the end of the timestep to provide us the corrected values:

$$u_i^{n+1} = u_i^n e^{-\frac{h}{\tau_i}} + \left(a_i \tau_i - s_i \tau_i^2\right)\left(1 - e^{-\frac{h}{\tau_i}}\right) + s_i \tau_i h. \tag{6}$$

With these 2a-b steps we obtain a second order method, which can be called "**linear-neighbour** method" and abbreviated by LN2.

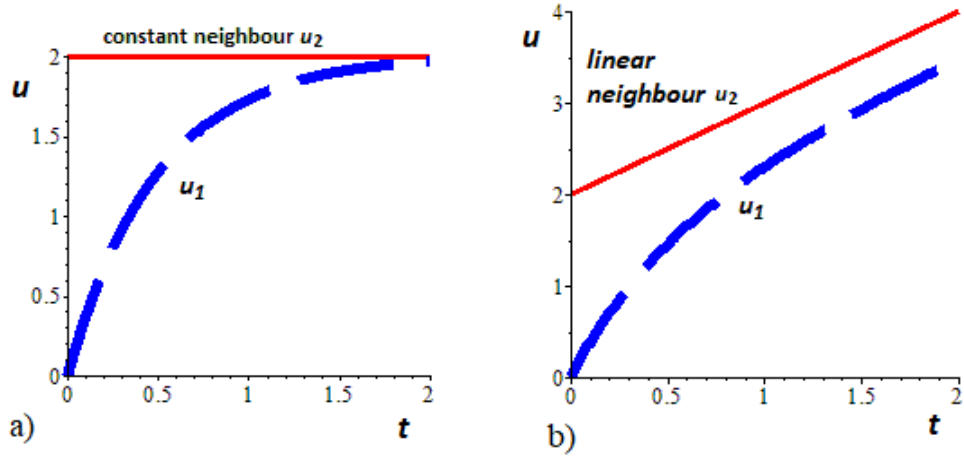

Figure 2 Visualization of the methods for a two-cell system. a) The continuous red line represents the $u_2$ variable in the second cell, which is considered constant during the timestep, while $u_1$, denoted by thick blue dashed line, is approaching it exponentially as in equation (4). b) At the second, linear-neighbour stage, $u_2$ is going up linearly to $u_2=4$ because of the external source, thus $u_1$ is approaching it according to formula (6).

It is possible to iterate both methods in the following way.

**3a)** Iteration of the constant-neighbour method. When, after completing step 1b), and we have $u_i^{n+1,\text{pred}}$, we calculate new $a_i$-s: $a_i = \sum_{j\in n(i)} m_{ij} u_j^{n+1,\text{pred}} + Q_i$. Now we can substitute these back to the ODE system (3) and after step 1b) again we obtain new $u_i^{n+1,\text{pred}}$ values. If we accepted it, this would be a two-stage CN method, briefly CN2. Then we can repeat this 3a) step again and again to get new $u_i^{n+1,\text{pred}}$ for CN3, CN4 etc., denoted by $u_i^{CN2}, u_i^{CN3},\ldots$ It means that the actual variable $u_i$ will, in each iteration step, exponentially tend to a new estimation of the neighbour's values.

**3b)** Iteration of the linear-neighbour method. After finishing step 2b) we can use the $u_i^{n+1}$ corrected values as predictors, obtain new $a_i^{\text{pred}}$ and $s_i$ and solve (5) again as it is explained at step 2b). If we do it once, we obtain a three-stage method (steps 1, 2 and 2 again), which can be denoted by LN3. We can go on with this iteration to obtain a four-stage method LN4, then LN5, etc. As the CN and LN methods have the same starting step 1a) and 1b), the LN1 method is the same as the CN1.

With these formulas we try to imitate the real physical processes in nature. In reality, if a system is thermally isolated, the temperature of each region of the system is approaching the equilibrium, which is the average temperature of the system. The speed of this process is inversely proportional to the heat capacity of the region and to the thermal resistance of the boundary of the region. We tried to apply this physical principle, even if it implies that the method cannot be made applicable generally for PDEs.

Regarding the above-mentioned convergence properties of the methods, we have the following convergence theorem, which is proved in Appendix A.

***Theorem 1.*** The CN1, CN2, …, (Constant Neighbour) methods are at least first order, the LN2, LN3, … (Linear Neighbour) methods are at least second order numerical methods for the general

$$\frac{d\vec{u}}{dt} = M\vec{u} + \vec{Q}, \quad \vec{u}(t=0) = \vec{u}(0)$$

linear ODE initial value problem, where *M* is an arbitrary constant matrix, while $\vec{Q}$ and $\vec{u}(0)$ are arbitrary vectors.

Before considering the stability properties of the schemes, we remind the reader that absolute stability means that arbitrarily large timesteps yield a bounded solution for the

$$\frac{du}{dt} = -\lambda u, \quad \lambda > 0$$

scalar equation [32, p. 13]. We note that no explicit Runge-Kutta methods fulfil this requirement which limits their use for solving PDEs [32, p. 24]. Now, for this scalar equation, i.e. in the case of one cell, our methods provide the analytical solution, thus, in this sense, they are trivially A-stable. We therefore examine the stability when $u$ is a vector, as it is always the case for spatially discretized PDEs. Unlike conventional explicit methods (like Runge-Kutta), formula (4) and (6) use a non-polynomial expression of the timestep-size $h$ to calculate the new values of the variable $u$. At the first stage, the temperature of each cell tends to the temperature of its neighbours, more precisely to $a_i \tau_i$, which is the weighted average of the temperatures of the neighbouring cells. This can be formulated in the following theorem.

**_Theorem 2._** *In case of the heat equation, when the $Q_i$ source terms are zero, the new values $u_i^{CN1}, u_i^{CN2}, u_i^{CN3}, \ldots$ and $u_i^{LN1}\left(=u_i^{CN1}\right), u_i^{LN2}, u_i^{LN3}, \ldots$ are the [convex combination]{.underline} of the initial values $u_j(0), \ j=1,\ldots,N$.*

This theorem is proved in Appendix B. That is why the result is always bounded by the initial extreme temperatures, which implies that the methods cannot be unstable for the heat equation, the stepsize can be arbitrarily large.

We summarize the most important properties of this class of methods in the following points.

1) They are obviously **explicit**; one can calculate the new values without solving an equation system or even without using matrices. It also implies that the process is easily parallelizable.

2) They are convergent, i.e. the solution converges to the exact solution of the discretized ODE system when the stepsize $h$ tends to zero. More specifically, the CN methods are at least **first order**, while the LN methods are at least **second order** in time.

3) We state that they are **stable** for heat conduction type problems, because the $u_i^{n+1}$ new value of the variable is the weighted average (convex combination) of $u_i^n$ and its neighbours $u_j^n$. Using these methods one can be sure that the solution automatically follows the Maximum and Minimum principles [32, p. 87], i.e. the extreme values of $u$ occur among the initial values (in the homogeneous $Q=0$ case, of course). This automatically implies that positivity is always preserved in the absence of heat sinks.

4) They are **one-step** methods in the sense that when we calculate $u_i^{n+1}$, we use only values $u_j^n$ at the beginning of the timestep, thus there is no need to store the $u_j^{n-1}$, $u_j^{n-2}$,… values of the variable at previous timesteps, i.e. the memory requirements are minimal.

5) They can be easily applied regardless of space dimensions, grid irregularity and inhomogeneity of the heat conduction medium.

6) According to our numerical experiments, the CN$i$ and LN$i$ iterations are always convergent when $i \to \infty$, however, they do not converge to the exact solution for any fixed timestep $h$. This statement is still waiting to be proved analytically; we managed to prove it only for special cases like a system consisting of two cells. We will show that a small number of iterations can significantly improve the accuracy, but no one should iterate too much in case of these methods.

We have to emphasize that these methods are **not** Finite Difference Methods (FDM) or Finite Volume Methods (FVM), even if they bear some resemblance to them. Here time derivatives are not approximated as in FDM and fluxes are not used as in FVM. Our methods are not conservative like FVM, and until now we were not able to find a conservative modification of it. Because of the exponential terms, our methods might also seem to be similar to the so-called exponential integrators. However, those methods use matrix exponentials [33], while we do not even need to use matrices during the calculation, thus our methods are fundamentally different from those.

### 3. Verification: comparison with an exact result

We solve PDE (1) on the interval $[0, \pi]$ with the following initial conditions:

$$u(x, t=0) = 10\sin(x) + 77\sin(2x)$$

and zero Dirichlet boundary conditions

$$u(x=0,t) = u(x=\pi,t) = 0$$

while $q \equiv 0$. It is easy to check that the analytical solution of this problem is

$$u(x,t) = 10\sin(x)\exp(-t) + 77\sin(2x)\exp(-4t)$$

We set the size of the spatial cells to $\Delta x = \pi/100 = 0.0314$, such that the first and the last cell-centre are at 0 and $\pi$ thus the number of cells are N=101. On Fig 3 one can see the global errors as a function of stepsize $h$. Error here means the maximum of the absolute value of the difference between the exact temperature and the temperature obtained by our numerical methods:

$$\max_i \left| u_i^{ex}(t_{fin}) - u_i^{num}(t_{fin}) \right|, \tag{7}$$

where, in this case, $t_{fin} = 1$.

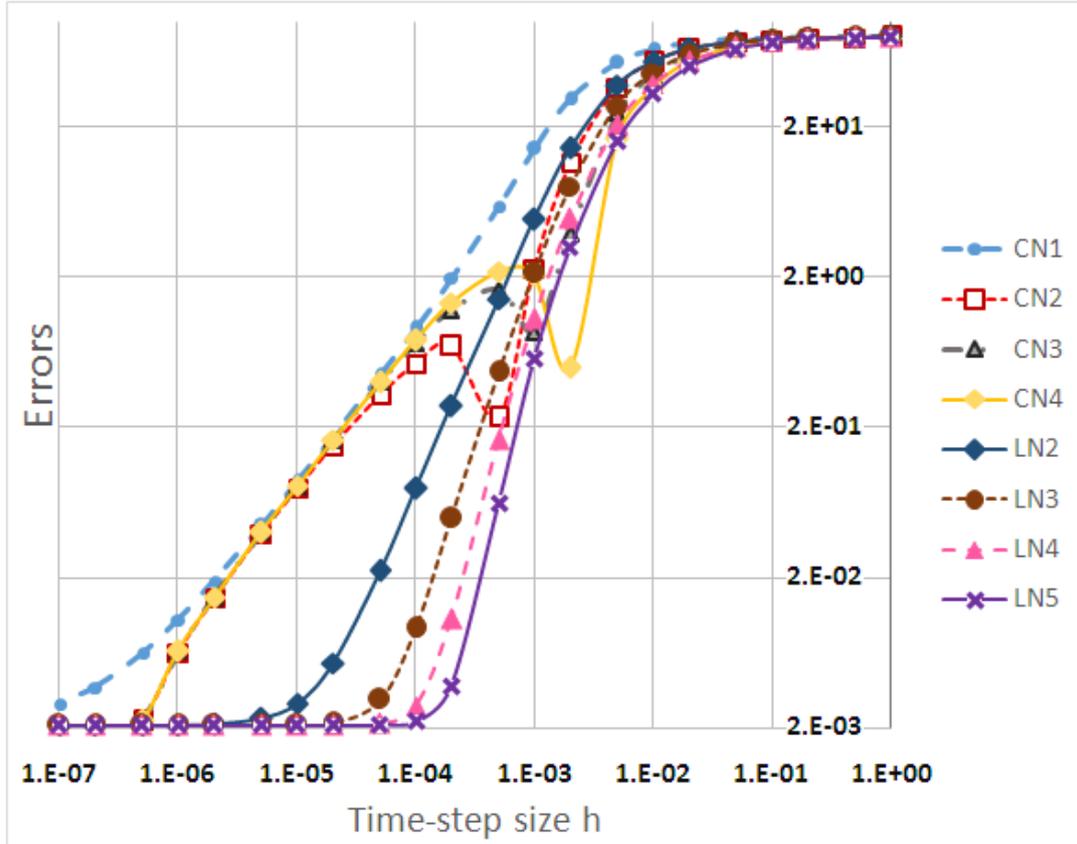

Figure 3 Maximum errors defined by formula (7) as the function of the timestep-size in case of our Constant- and Linear-neighbour methods.

One can see that the error tends to a small nonzero value in all cases. This non-vanishing error is caused by the discretization of the space variable and tends to zero with $\Delta x \to 0$. We have to note here that the thorough investigation of the consistency of the methods must be done in the future. From the figure one might think that the LN3-5 methods have an order larger than 2. However, later we will see that this is not the case. We remind the reader that because of the exponential term, the expressions of the new variables $u_i^{n+1}$ contain the stepsize $h$ up to infinite order and that can cause apparently different orders.

## 4. Comparison with numerical results 1

First we examined a regular rectangle-shaped lattice with $50 \times 20 = 1000$ cells. Different random values were given to the capacities and to the resistances:

$$C_i = 10^{(1-2 \cdot rand)}, \quad R_{x,i} = 10^{(1-2 \cdot rand)}, \quad R_{y,i} = 10^{(1-2 \cdot rand)},$$

where *rand* is a random number generated by the MATLAB uniformly in the (0,1) interval for each quantity. It means that these quantities followed a log-uniform distribution between 0.1 and 10. The initial temperatures had a uniform distribution between 0 and 1000, $u_i(0) = 1000 \cdot (1 - rand)$ while the source-terms had a uniform distribution between -500 and 500, $Q_i = 1000 \cdot (rand - 0.5)$. The task was to solve this system for the temperatures between $t_0$=0s and $t_{FIN}$=1s.

The stiffness ratio, i.e. the ratio of the largest and smallest (nonzero) eigenvalues of the matrix, is $5.36 \cdot 10^5$. For the explicit Euler method (which is equivalent to the forward-time central-space FTCS scheme for the original PDE), the maximum possible timestep is

$$h_{MAX}^{EE} = \left|\frac{2}{\lambda_m}\right| = 0.0085,$$

above this threshold instability necessarily occurs. Here $\lambda_m$ is (non-positive) eigenvalue of the matrix with the largest absolute-value. Albeit this problem is not extremely stiff, in order to obtain a reference solution, we used an implicit ode15s solver of MATLAB with strict error tolerance ('RelTol' and 'Abstol' were $10^{-12}$), where the letter *s* indicates that the codes were designed especially for stiff systems. In Fig. 4 and 5 we present the maximum error for the Constant-neighbour and the Linear-neighbour methods in case of different number of iterations as a function of the stepsize *h* while in Fig. 6 one can see all curves for large values of *h*. From the figures one can see that for small h, the global error decreases with the first power of the stepsize in case of the Constant-neighbour methods and with the second power for the Linear-neighbour methods, which confirms that the methods are first and second order respectively, except of course LN1, which is the same as CN1.

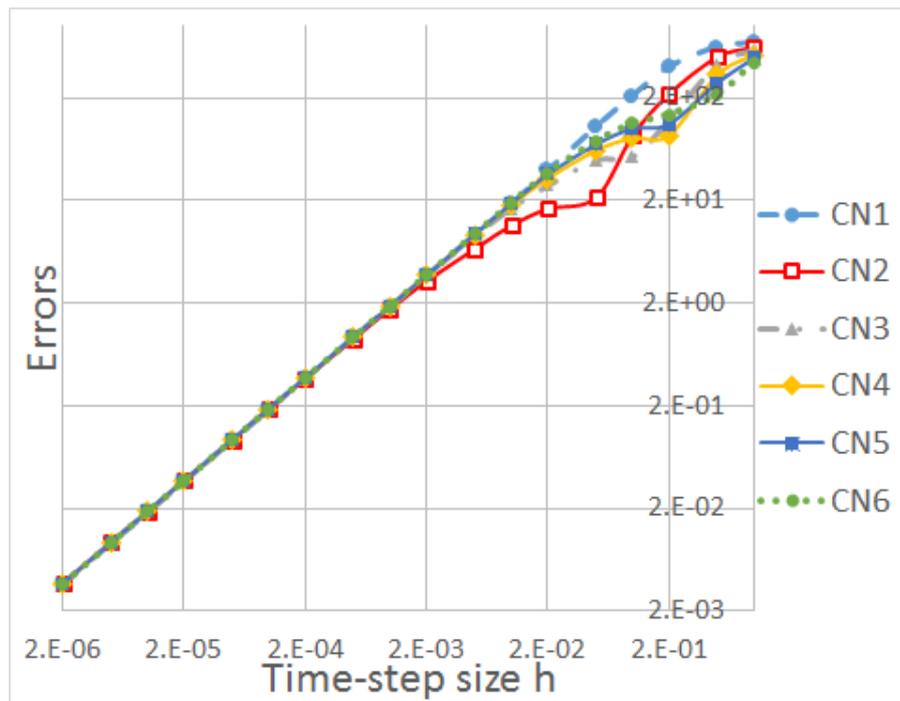

Figure 4 Maximum errors defined by formula (7) as the function of the timestep-size *h* in case of our Constant-neighbour methods. For *h*<0.001 the curves are almost identical, thus in this region the iteration makes no sense.

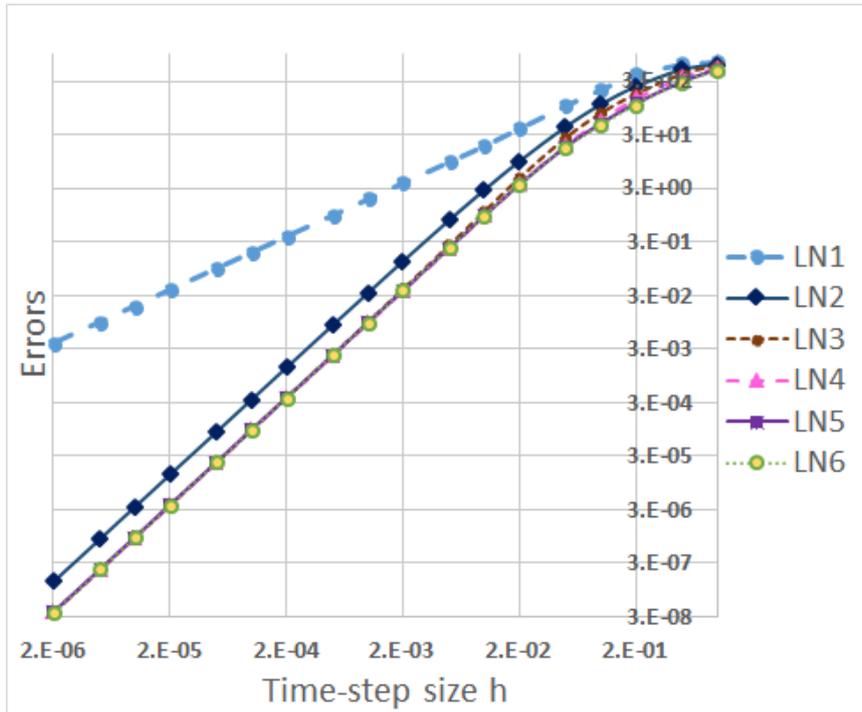

Figure 5 Maximum errors as the function of the timestep-size in case of our Linear-neighbour methods. We remind the reader that LN1 is the same as CN1.

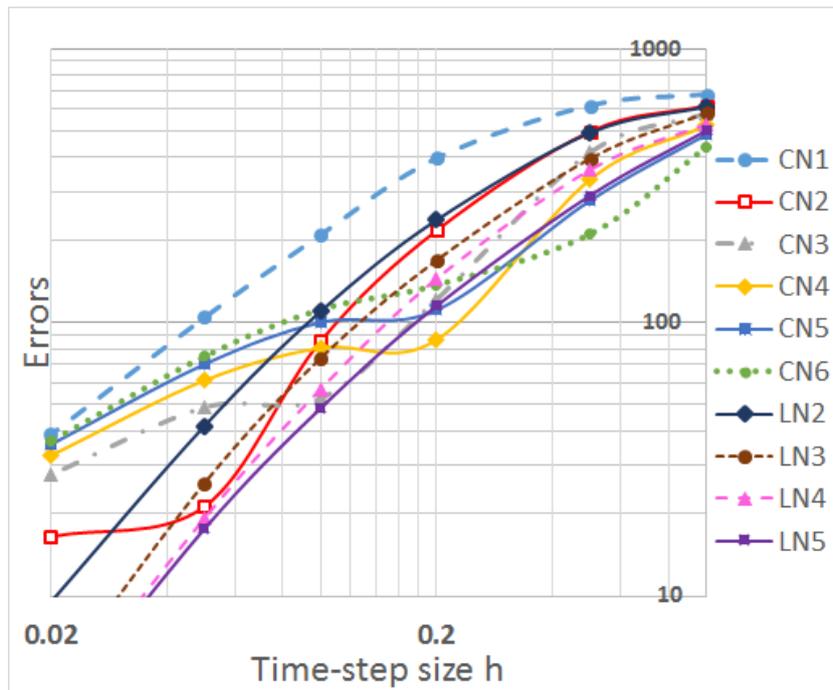

Figure 6 Maximum errors as the function of the timestep-size in case of our Constant- and Linear-neighbour methods with different number of iterations for larger timestep-sizes.

We also examined the errors as the function of the number of iterations for fixed timestep-sizes. The results are presented in Fig. 7. It is interesting that the decrease of errors is not monotonous in the case of CN methods. We admit that currently we don't know the reason of this anomalous behaviour. At some specific cases a CN method performs extraordinary well, but it seems to be impossible to predict under which circumstances we can expect this.

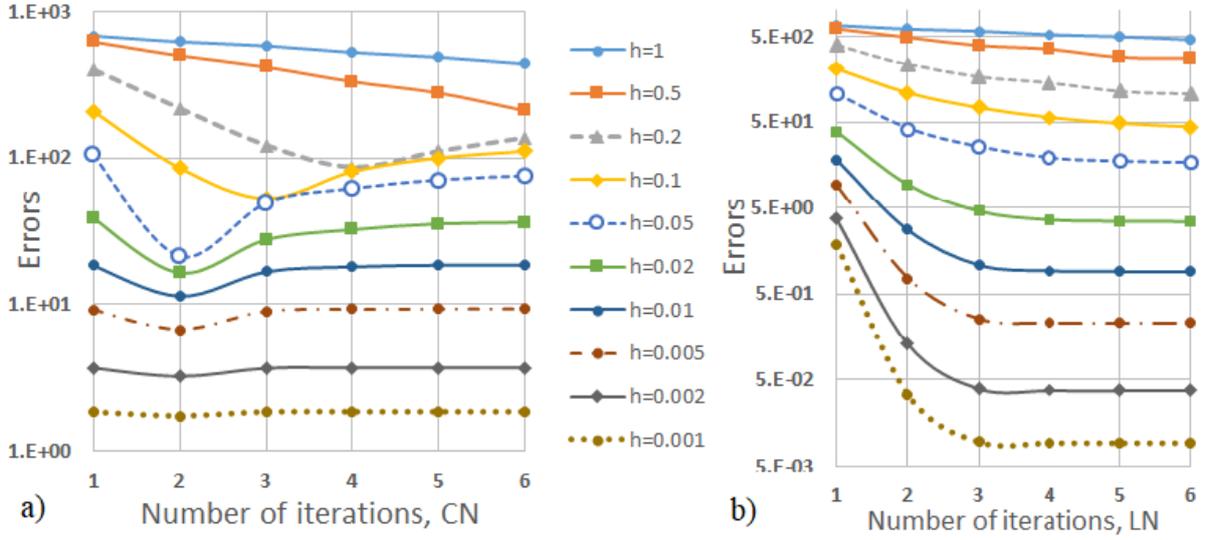

Figure 7 Maximum errors of CN and LN methods as the function of iteration number for different timestep-size

On the other hand, we can conclude that in case of Linear Neighbour methods, iteration systematically improves the precision. However, we still have to answer the question that more iterations or decreasing the timestep is a better decision. To answer that, we fixed the total number of iterations GITC which is the product of the number of iterations (or stages) per timestep and the number of timesteps. This quantity is almost strictly proportional to the total runtime, and we can examine the errors as a function of the number of iterations per timestep. In order to obtain more balanced information about the errors, we introduce two new error quantities. The first one, $SumD$ is the sum of the deviations for all of the cells:

$$SumD = \sum_{i=1}^{N} \cdot \left| u_i^{ex}(t_{fin}) - u_i^{num}(t_{fin}) \right| \tag{8}$$

The second one is the sum of the same differences but weighted with the capacities,

$$SEnD = \sum_{i=1}^{N} C_i \cdot \left| u_i^{ex}(t_{fin}) - u_i^{num}(t_{fin}) \right|, \tag{9}$$

and thus it gives the error in terms of energy. It reflects the fact that a temperature deviance in a big cell has more significance in practice than in a tiny cell. As the number of cells is large, these errors are usually much larger than the previously defined maximum error. To make the graphical representation of them easier, we normalized them by the square root of the number of cells:

$$SumDN = \frac{SumD}{\sqrt{N}}, \quad SEnDN = \frac{SEnD}{\sqrt{N}} \tag{10}$$

First, the Global Iteration Counter GITC were fixed to 420, the least common multiple of the first 7 positive integer numbers. It means that for one iteration-step (CN1=LN1) the timestep-size was $t_{fin}/420 = 0.00238$, for LN2 it was twice as large, etc. The results can be seen in Fig. 8.

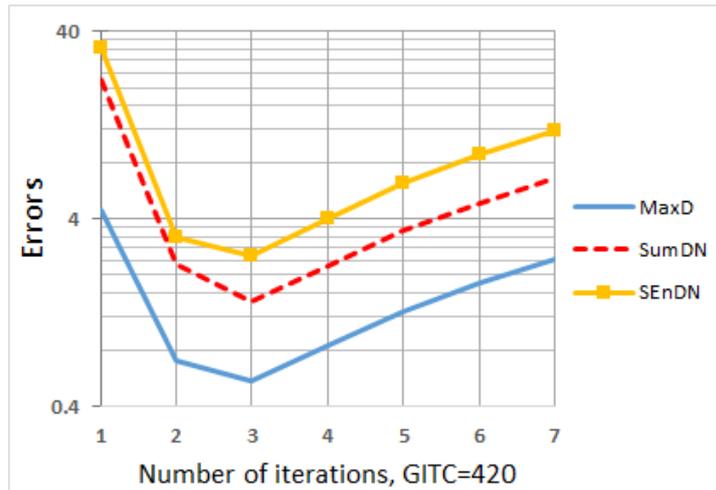

Figure 8 Different errors of LN methods defined by formulas (7) and (10) as the function of the iteration number for fixed Global Iteration Counter (GITC=420). We emphasize that the timestep-size is also increasing with the number of iterations (from left to right in the diagram), thus the total running time is nearly constant.

We repeated this for larger global number of iterations (GITC=2100 and 10500) and presented the maximum errors in Fig. 9.

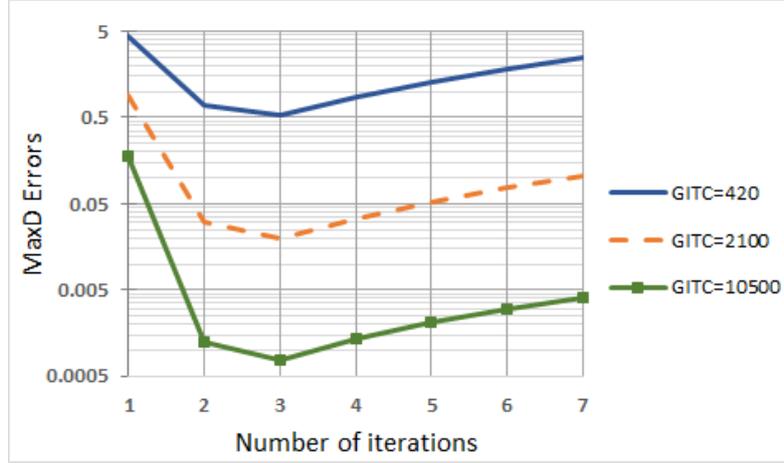

Figure 9 Maximum errors of LN methods as the function of the iteration number for 3 different Global Iteration Counter value.

Looking at Figures 8 and 9 we can conclude that it is worth to iterate three (one constant and then two linear-neighbour stages) but not more, i.e. the LN3 method has the best performance.

In Table 1 we summarize some results obtained by MATLAB routines ode15s, ode23s, ode23, ode23t (with very loose error tolerance to increase the speed) and our methods CN2 and LN3. One can see that our methods are much faster than the conventional explicit or implicit methods, even without any optimization, adaptive stepsize control, parallelization or vectorization.

## 5. Comparison with numerical results 2

The second system was a rectangle-shaped lattice with $250 \times 20 = 5000$ cells. The capacities and the resistances followed a log-uniform distribution between 0.001 and 1000:

$$C_i = 10^{(3-6 \cdot rand)}, \ R_{x,i} = 10^{(3-6 \cdot rand)}, \ R_{y,i} = 10^{(3-6 \cdot rand)},$$

The initial temperatures were zero everywhere, $u_i(0) = 0$ while the source-terms had a uniform distribution between -500 and 500, $Q_i = 1000 \cdot (rand - 0.5)$. The task was to solve this system for the temperatures between $t_0$=0s and $t_{FIN}$=1s.

The stiffness ratio is rather high, $1.12 \cdot 10^{12}$. For the explicit Euler method the maximum possible timestep is $h_{MAX}^{EE} = 1.1 \cdot 10^{-6}$. To obtain a reference solution, we used the same solver with the same error tolerance as in the previous section. In Fig. 10 and 11 we present the maximum error for the Constant-neighbour and the Linear-neighbour methods in case of different number of iterations as a function of the stepsize $h$.

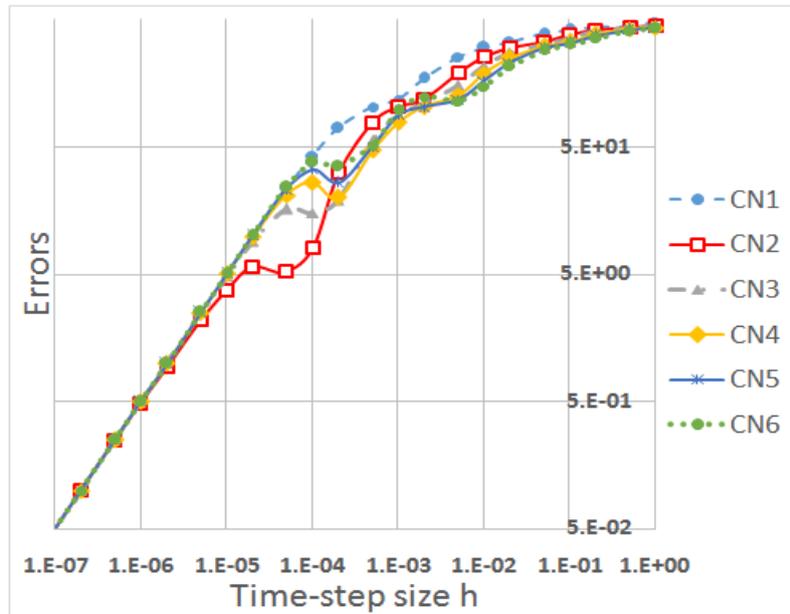

Figure 10 Maximum errors defined by formula (7) as the function of the timestep-size in case of the Constant-neighbour methods.

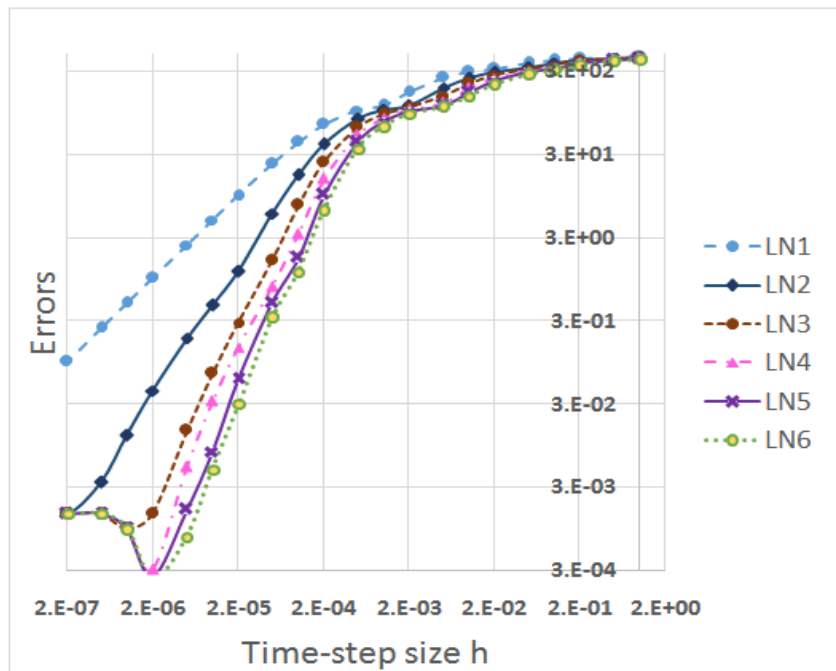

Figure 11 Maximum errors as the function of the timestep-size in case of the Linear-neighbour methods. We remind the reader that LN1 is the same as CN1.

One can see that the anomaly which was present only in the case of the CN methods in the previous section are now stronger and something similar appeared for the LN methods as well. To reveal the nature of these anomalies, we plan to test the methods for much larger and stiffer systems as well, but this will be possible only by using more powerful computers than we currently have.

We also examined the errors as the function of the number of iterations for fixed timestep-sizes. The results are presented in Fig. 12.

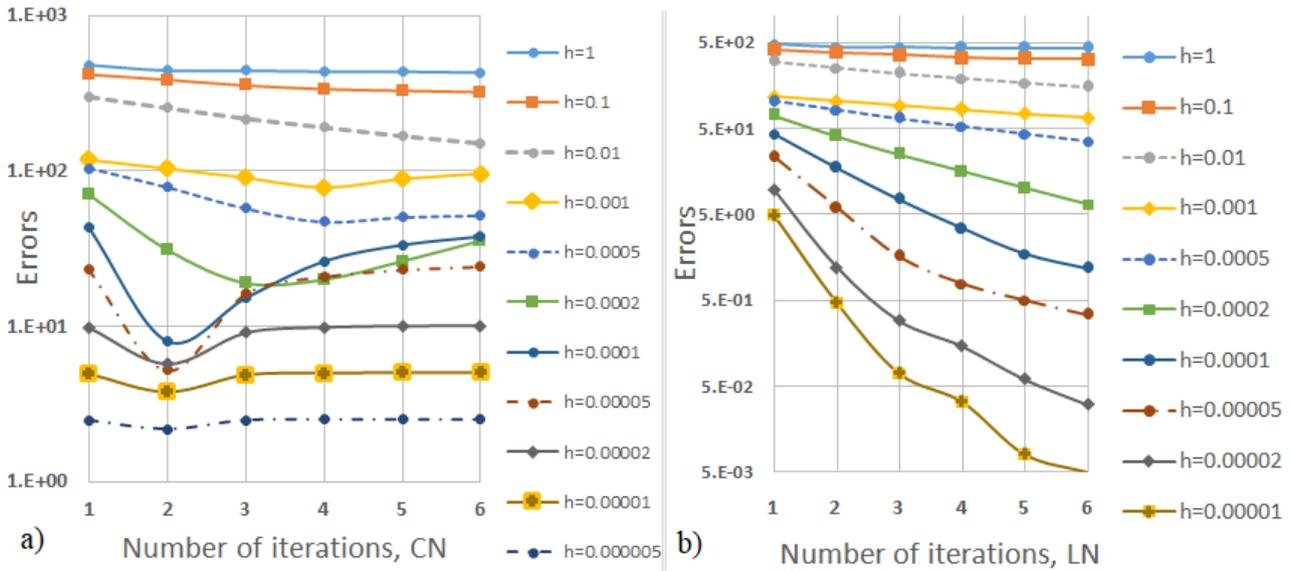

Figure 12 Maximum errors of CN methods as the function of the iteration number for different timestep-size

As in the previous case, one can conclude that (unlike in the case of the CN method,) the iteration of the Linear Neighbour method normally improves the precision. We also examined the errors for fixed total number of iterations GITC, which first was set to 840, then multiplied subsequently by five. The results can be seen in Fig. 13.

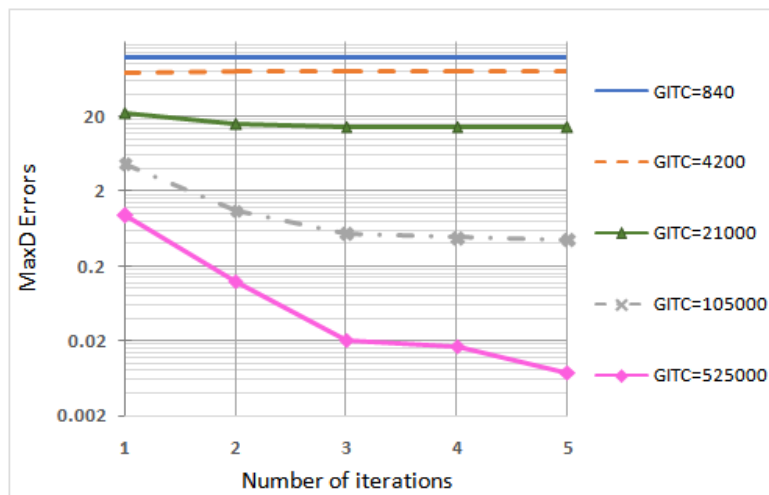

Figure 13 Maximum errors of LN methods as the function of the iteration number for 3 different Global Iteration Counter

In Table 2. we summarize some results obtained by MATLAB routines ode15s, ode23s and ode23t (with very loose error tolerance) and our methods CN2 and LN3. We note that the explicit solvers ode23 and ode45 were not able to provide any results, because for large error tolerance they failed to converge, while tightening the tolerance resulted in an error message "array exceeds maximum array size" after hours of running.

One can see that our methods are much faster than the conventional explicit or implicit methods if we set the same moderate precision requirements. One might think that the method CN2 for $h$=0.002 produces too large errors. However, the average error is only 2.97 and from Fig. 14 one can see that the results fit well to the reference curve.

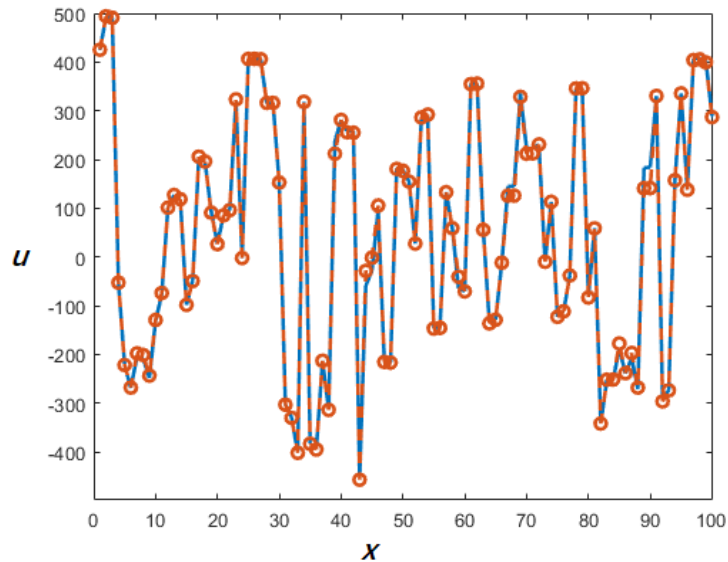

Figure 14 The variable *u* as a function of the space variable *x* for the first 100 cells. The blue line represents the high-precision solution while the orange dashed line with circles are the values produced by our CN2 algorithm for $h=0.002$ with orders of magnitude faster than the conventional solvers.

Bearing in mind that our methods can produce acceptable results with orders of magnitude faster than the well-established MATLAB routines, we emphasize that this is still a small system with 5000 cells, and for larger systems implicit methods have more serious drawbacks. We note that the reason we have not used larger number of cells for testing purposes is that our computer cannot solve them by implicit methods because of the too large memory requirements.

## 6. Conclusions and future work

In this article we proposed a class of new explicit numerical algorithms to solve the spatially discretized non-stationary heat conduction or diffusion equation with external sources. The predictor "constant-neighbour" method and its iterations are first order in time, while the "linear-neighbour" corrector step and its iterations provide second order methods. In the absence of external sources, the variables at the end of the timestep are the weighted average (convex combinations) of the values at the beginning of the timestep, thus the method is stable for arbitrarily large timestep-sizes. We presented the exact analytical proof of these statements in the Appendices.

We verified the algorithms by a comparison to a simple analytical solution of the original PDE, and then illustrated the performance of the methods in the case of two discretized 2D systems with inhomogeneous random parameters and initial conditions. According to numerical experiments, the proposed schemes are competitive, even without any optimization, adaptive stepsize control, parallelization or vectorization. Increasing the number of cells gives us increasing advantage over the implicit schemes while increasing stiffness makes our methods more competitive against the (conditionally stable) explicit schemes. However, we do not state that these methods are ideal when someone needs very high precision results. We usually found that the LN methods are significantly more accurate, but they are still only second order methods. We have already begun to seek higher-order versions of the methods. Besides this we have started or plan to start working in the following directions:

- Systematically investigate the consistency of the methods.
- Test the methods in case of real-life problems.
- Examine how parallel programming can accelerate the calculation.
- Develop the adaptive stepsize control version of the algorithms.
- Elaborate how they can be combined with the geometric multigrid methods.
- Extend the methods to more complicated equations like the advection-diffusion equation or to nonlinear cases when the parameters like the specific heat depends on the temperature as well.

To make progress towards these goals and implement other ideas, we are open to scientific collaboration.

## 7. Acknowledgement

This research was supported by the EU and the Hungarian State, co-financed by the ERDF in the framework of the GINOP-2.3.4-15-2016-00004 project. I'm also grateful to András Gilicz, László Baranyi, Gábor Pszota and Tamás Jenyó.

| method | Runtime (s) | MaxD | SumD | SEnD |
|---|---|---|---|---|
| ode45 | 0.55 | 25.8 | 30 | 5.78 |
| ode15s | 1.56 | 9.26 | 2078.9 | 4688 |
| ode23s | 14.54 | 0.31 | 58.2 | 109.3 |
| ode23 | 0.3 | >$10^5$ (unstable) | | |
| ode23 | 0.4 | 0.025 | 0.13 | 0.15 |
| ode23t | 1.89 | 0.45 | 84.1 | 157.9 |
| CN2, $h=0.02$ | 0.006 | 16.5 | 2966 | 4933 |
| LN3, $h=0.05$ | 0.027 | 25.7 | 2239 | 3575 |
| LN3, $h=0.01$ | 0.068 | 1.08 | 91.8 | 158.8 |

Table 1. Performance and errors of different solvers for N=1000 cells. CN is for "Constant-Neighbour", LN is for "Linear-Neighbour". The error quantities are defined by formulas (7-9).

| method | Runtime (s) | MaxD | SumD | SEnD |
|---|---|---|---|---|
| ode15s | 171.4 | 8.17 | 4788 | 172112 |
| ode23s | 1816 | 0.187 | 72.03 | 1960 |
| ode23t | 163.1 | 0.377 | 153.1 | 3905 |
| CN2, $h=0.002$ | 0.63 | 117.8 | 14863 | 95649 |
| CN1, $h=0.0002$ | 1.75 | 70.48 | 4380 | 27110 |
| CN2, $h=0.0001$ | 6.53 | 8.0 | 629.4 | 4911 |
| LN4, $h=0.0001$ | 12.5 | 3.4 | 102.9 | 326 |
| LN3, $h=0.00002$ | 49.6 | 0.292 | 11.38 | 45.3 |
| LN3, $h=0.00001$ | 95.1 | 0.072 | 2.78 | 12.6 |

Table 2. Performance and errors of different solvers for N=5000 cells. CN is for "Constant-Neighbour", LN is for "Linear-Neighbour". The error quantities are defined by formulas (7-9).

# 9. Appendix A. The proof of Theorem 1

To help the reader follow the argument, we present the matrix form of the ODE system:

$$\frac{d}{dt}\begin{pmatrix} u_1 \\ \vdots \\ u_N \end{pmatrix} = \begin{pmatrix} m_{11} & \cdots & m_{1N} \\ \vdots & \vdots & \vdots \\ m_{N1} & \cdots & m_{NN} \end{pmatrix} \begin{pmatrix} u_1 \\ \vdots \\ u_N \end{pmatrix} + \begin{pmatrix} Q_1 \\ \vdots \\ Q_N \end{pmatrix}$$

The exact solution of the initial value problem:

$$\vec{u}(t) = e^{Mt}\vec{u}_0 + \left(e^{Mt} - 1\right)M^{-1}\vec{Q} = \left(1 + Mt + M^2\frac{t^2}{2} + M^3\frac{t^3}{3!} + \ldots\right)\vec{u}_0 + \left(t + M\frac{t^2}{2} + M^2\frac{t^3}{3!} + \ldots\right)\vec{Q}$$

Without the loss of generality we examine only $u_1$ at the first timestep. It implies that we use the initial values $u_j(0)$ at the beginning of the time-step. The 0th and first order terms in the exact solution at $t=h$ are the following:

$$u_1(h) = u_1(0)\left(1 + m_{11}h\right) + h\sum_{j>1} m_{1j}u_j(0) + Q_1 h \tag{A1}$$

The second order terms:

$$u_1(h) = \frac{h^2}{2}\sum_{j=1}^{N} m_{1j}\sum_{k=1}^{N} m_{jk}u_k(0) + \frac{h^2}{2}\sum_{j=1}^{N} m_{1j}Q_j \tag{A2}$$

A) Our **constant**-neighbour solution:

$$u_1^{CN}(h) = u_1(0)e^{m_{11}h} + \frac{a_1}{m_{11}}\left(e^{m_{11}h} - 1\right),$$

where

$$a_1 = \sum_{j>1} m_{1j}u_j(0) + Q_1$$

is a constant. Now one can write

$$u_1^{CN}(h) = u_1(0)\left(1 + m_{11}h + \frac{m_{11}^2 h^2}{2} + \ldots\right) + \frac{1}{m_{11}}\left(\sum_{j>1} m_{1j}u_j(0) + Q_1\right)\left(m_{11}h + \frac{(m_{11}h)^2}{2} + \ldots\right) =$$

$$u_1(0)\left(1 + m_{11}h + \frac{m_{11}^2 h^2}{2} + \ldots\right) + h\left(\sum_{j>1} m_{1j}u_j(0) + Q_1\right)\left(1 + \frac{m_{11}h}{2} + \ldots\right) \tag{A3}$$

It is easy to see that the first and the second order terms are the same as in (A1), thus the local error is second order, which means that the CN method is first order.

B) When we do the iteration in the case of the CN method, we put the newly obtained values $u_j^{CN}(h)$ instead of the initial values $u_j(0)$ in (A3). However, the term with $u_j(0)$ is already multiplied by $h$ in (A3), thus we need to care only about the 0th order terms in $u_j^{CN}(h)$, which is just $u_j(0)$. It means that the iterations modify only the terms with order higher than one, so the iterated versions CN2, CN3, etc. are at least first order.

C) In the case of the second (**linear**-neighbour) method, we start from (6)

$$u_1^{LN2} = u_1(0)e^{-\frac{h}{\tau_1}} + \left(a_1\tau_1 - s_1\tau_1^2\right)\left(1 - e^{-\frac{h}{\tau_1}}\right) + s_1\tau_1 h = u_1(0)e^{m_{11}h} + \frac{a_1}{m_{11}}\left(e^{m_{11}h} - 1\right) + \frac{s_1 h}{m_{11}}\left(\frac{e^{m_{11}h} - 1}{m_{11}h} - 1\right).$$

The first two terms are obviously the same as $u_1^{CN}(h)$ in point A). After calculating the quantities

$$a_1^{pred} = \sum_{j>1} m_{1j}u_j^{CN}(h) + Q_1 \quad \text{and} \quad s_1 = \frac{a_1^{pred} - a_1}{h} = \frac{1}{h}\sum_{j>1} m_{1j}\left(u_j^{CN}(h) - u_j(0)\right)$$

we can examine the last term:

$$\frac{s_i h}{m_{11}}\left(\frac{e^{m_{11}h}-1}{m_{11}h}-1\right) = \frac{1}{m_{11}}\sum_{j>1} m_{1j}\left(u_j^{CN}(h)-u_j(0)\right)\cdot\left(\frac{m_{11}h+\frac{m_{11}^2 h^2}{2}+\frac{m_{11}^3 h^3}{6}+\ldots}{m_{11}h}-1\right) = \qquad (A4)$$

$$= \sum_{j>1} m_{1j}\underbrace{\left(u_j^{CN}(h)-u_j(0)\right)}_{X}\cdot\underbrace{\left(\frac{h}{2}+\frac{m_{11}h^2}{6}+\ldots\right)}_{*}$$

Now we substitute back the $u_j^{CN}(h)$ solution instead of (X) into the first bracket:

$$u_j^{CN}(h)-u_j(0) = u_j(0)\left(1+m_{jj}h+\ldots\right)+h\left(\sum_{k\neq j} m_{jk}u_k(0)+Q_j\right)\left(1+\frac{m_{jj}h}{2}+\ldots\right)-u_j(0) =$$

$$u_j(0)\left(m_{jj}h+\ldots\right)+h\left(\sum_{k\neq j} m_{jk}u_k(0)+Q_j\right)\underbrace{\left(1+\frac{m_{jj}h}{2}+\ldots\right)}_{\square}$$

This term contains only first and higher order terms, as well as the last (*) bracket in (A4). We can conclude that the LN2 iteration step does affect only the second and higher order terms. As we do not examine third and higher order terms, we have to care only for the first, $h/2$ term in (*) and also the first term, the number one in ($\square$), thus we can write for the second order terms of (A4):

$$\sum_{j>1} m_{1j}\left[u_j(0)\left(m_{jj}h+\ldots\right)+h\left(\sum_{k\neq j} m_{jk}u_k(0)+Q_j\right)\right]\cdot\frac{h}{2} = \frac{h^2}{2}\sum_{j>1} m_{1j}\left[u_j(0)m_{jj}+\sum_{k\neq j} m_{jk}u_k(0)+Q_j\right]$$

The second order terms from (A3) and (A4) are the following:

$$u_1(0)\frac{m_{11}^2 h^2}{2}+h\left(\sum_{j>1} m_{1j}u_j(0)+Q_1\right)\frac{m_{11}h}{2}+\frac{h^2}{2}\left[\sum_{j>1} m_{1j}m_{jj}u_j(0)+\sum_{j>1} m_{1j}\sum_{k\neq j} m_{jk}u_k(0)+\sum_{j>1} m_{1j}Q_j\right] =$$

$$\frac{h^2}{2}\Bigg[\underbrace{u_1(0)m_{11}^2+m_{11}\sum_{j>1} m_{1j}u_j(0)}_{m_{11}\sum_{j=1}^{N} m_{1j}u_j(0)}+\underbrace{\sum_{j>1} m_{1j}m_{jj}u_j(0)+\sum_{j>1} m_{1j}\sum_{k\neq j} m_{jk}u_k(0)}_{\sum_{j>1} m_{1j}\sum_{k=1}^{N} m_{jk}u_k(0)}+\underbrace{m_{11}Q_1+\sum_{j>1} m_{1j}Q_j}_{\sum_{j=1}^{N} m_{1j}Q_j}\Bigg],$$

$$\underbrace{\hspace{6cm}}_{\sum_{j=1}^{N} m_{1j}\sum_{k=1}^{N} m_{jk}u_k(0)}$$

which is equivalent to (A2), i.e. the LN2 method is at least 2$^{nd}$ order.

D) We can follow the same logic as in point B). When we calculate the value of $u_1^{LN3}$, we have to replace $u_j^{CN}(h)$ with $u_j^{LN2}$ in $a_1^{pred}$ and therefore in the (X) term in (A4). Again, these $u_j^{LN2}$ are multiplied by (*) containing at least first order terms, thus we have to take only 0$^{th}$ and 1$^{st}$ order terms into consideration in $u_j^{LN2}$. As the LN2 iteration step does not change the 0$^{th}$ and 1$^{st}$ order terms, we obtain that the LN3 (and also the LN4, etc.) results equal to the LN2 results up to 2$^{nd}$ order, which means that we proved the theorem.

## 10. Appendix B. The proof of Theorem 2

Again, without the loss of generality we examine only $u_1$ at the first timestep. We will use the facts that $\tau_i = -1/m_{ii} > 0$, $m_{ij,j\neq i} > 0$ and therefore $0 < e^{-h/\tau_i} = e^{m_{ii}h} \leq 1$ because of physical reasons, namely the Second law of thermodynamics.

A) Let us examine the CN1 solution:

$$u_1^{CN} = u_1(0)\cdot e^{-\frac{h}{\tau_1}}+a_1\tau_1\cdot\left(1-e^{-\frac{h}{\tau_1}}\right)$$

The term $a_1 \tau_1$ is the following:

$$a_1 \tau_1 = \sum_{j>1} m_{1j} u_j^n \cdot \frac{-1}{m_{11}} = \sum_{j>1} \frac{u_j(0)}{C_1 R_{1j}} \frac{C_1}{\sum_{k>1} \frac{1}{R_{1k}}} = \frac{\sum_{j>1} \frac{u_j(0)}{R_{1j}}}{\sum_{k>1} \frac{1}{R_{1k}}}$$

thus it is a convex combination of initial values $u_j(0)$, and therefore $u_1^{CN}$ is also a convex combination of the initial values $u_j(0)$, $j=1,...,N$. Indeed, the coefficients $e^{-\frac{h}{\tau_1}}$ and $\frac{1/R_{1j}}{\sum_{k>1} 1/R_{1k}} \left(1 - e^{-\frac{h}{\tau_1}}\right)$ are nonnegative and their sum is one.

B) When we perform the iteration in the case of the CN method, we put $u_j^{CN}(h)$ instead of the initial value $u_j(0)$. However, as we have just proven, the term $u_j^{CN}(h)$ is a convex combination of $u_k(0)$, thus $u_j^{CN2}(h)$ will be a convex combination as well, which can be said also about the CN3, CN4 etc. results.

C) At the case of the LN2 (**linear-neighbour**) method, we start from the formulas we had already obtained in Appendix A, point C).

$$u_1^{LN2} = u_1(0) e^{m_{11} h} + \frac{a_1}{m_{11}} \left(e^{m_{11} h} - 1\right) + \frac{s_1 h}{m_{11}} \left(\frac{e^{m_{11} h} - 1}{m_{11} h} - 1\right), \tag{B1}$$

where

$$s_1 = \frac{1}{h} \sum_{j>1} m_{1j} \left(u_j^{CN}(h) - u_j(0)\right) \tag{B2}$$

The first two terms in (B1) are the same as in the case of the CN1 method, thus we know that they are the convex combination of the initial values $u_j(0)$. Now let us examine the bracket in (B2). The term $u_j^{CN}(h)$ is a convex combination of the initial values $u_k(0)$, thus the sum of the coefficients is 1. The coefficient of the second term $u_j(0)$ is $-1$, thus the bracket is a combination of the initial values where the sum of the coefficients is zero. The remaining question is whether the coefficient of any concrete $u_j(0)$ in the (B1) expression of $u_1^{LN2}$ is nonnegative. It is enough to examine the second and the third term:

$$\frac{1}{m_{11}} \sum_{j>1} m_{1j} u_j(0) \left(e^{m_{11} h} - 1\right) + \sum_{j>1} m_{1j} \left(u_j^{CN}(h) - u_j(0)\right) \frac{1}{m_{11}} \left(\frac{e^{m_{11} h} - 1}{m_{11} h} - 1\right).$$

Substituting the formula for $u_j^{CN}(h)$ we can write the coefficient of $u_j(0)$ as

$$\frac{m_{1j}}{-m_{11}} \left[\left(1 - e^{m_{11} h}\right) + \left(1 - e^{m_{jj} h}\right) \left(\frac{e^{m_{11} h} - 1}{m_{11} h} - 1\right)\right]$$

Now the coefficient in front of the square bracket is nonnegative while the expression in the square bracket is also nonnegative for $m_{ii} < 0$, $m_{jj} \leq 0$ and $h > 0$, which means that $u_1^{LN2}$ is again a convex combination of the initial values.

D) When examining the LN3 and higher methods, we can combine the arguments of the previous points. When we calculate the value of $u_1^{LN3}$, we have to replace $u_j^{CN}(h)$ with $u_j^{LN2}$ in (B2) and therefore (B1), which are already convex combinations of the initial values, therefore the LN3, LN4, etc. iterations leave the convex combination property unchanged.